\pdfoutput=1
\documentclass[12pt]{article}

\usepackage[T1]{fontenc}

\usepackage[utf8]{inputenc}

\usepackage[english]{babel}

\usepackage[nomath]{lmodern}

\usepackage{amsmath}
\usepackage{amssymb}
\usepackage{amsthm}

\usepackage{stmaryrd}

\usepackage{geometry}

\usepackage{enumitem}

\usepackage{authblk}

\usepackage{tikz}
\usetikzlibrary{arrows.meta} 
\usetikzlibrary{graphs} 
\usetikzlibrary{decorations} 

\theoremstyle{definition}

\newtheorem{theorem}{Theorem}[section]
\newtheorem{definition}[theorem]{Definition}

\newtheorem{proposition}[theorem]{Proposition}
\newtheorem{observation}[theorem]{Observation}
\newtheorem{corollary}[theorem]{Corollary}

\newtheorem{remark}[theorem]{Remark}
\newtheorem{lemma}[theorem]{Lemma}
\newtheorem{example}[theorem]{Example}

\newtheorem{question}[theorem]{Question}
\newtheorem*{acknowledgements}{Acknowledgements}

\newcommand*\cdef{\newcommand*}

\setlist{itemsep = 0pt}
\setlist[enumerate]{leftmargin=*} 
\setlist[enumerate, 1]{label=\upshape (\roman*), ref=(\roman*)}
\setlist[enumerate, 2]{label=\upshape (\alph*), ref=(\alph*)}

\cdef \DescriptionFormat [1]{%
	\normalfont\emph{#1}%
}

\setlist[description]{format=\DescriptionFormat}

\cdef \Show [1]{%
	\expandafter\show\csname #1\endcsname
}

\cdef \AlignedMath [1]{%
	\begingroup
	
	\begin{align*}#1\end{align*}%
	\endgroup
	\ignorespaces
}

\let \O \undefined 
\let \P \undefined 
\let \H \undefined 

\begin{document}

\cdef \DefUnicode [2]{%
	\expandafter\cdef
		\csname u8:\detokenize{#1}\endcsname
		{#2}%
}

\cdef \ShowUnicode [1]{%
	\expandafter\show
		\csname u8:\detokenize{#1}\endcsname
}

\DefUnicode{¬}{\lnot}

\DefUnicode{×}{\times}


\DefUnicode{α}{\alpha}
\DefUnicode{β}{\beta}
\DefUnicode{γ}{\gamma}
\DefUnicode{δ}{\delta}
\DefUnicode{ε}{\varepsilon}
\DefUnicode{ζ}{\zeta}
\DefUnicode{η}{\eta}
\DefUnicode{θ}{\theta}
\DefUnicode{ι}{\iota}
\DefUnicode{κ}{\kappa}
\DefUnicode{λ}{\lambda}
\DefUnicode{μ}{\mu}
\DefUnicode{ν}{\nu}
\DefUnicode{ξ}{\xi}
\DefUnicode{π}{\pi}
\DefUnicode{ρ}{\rho}
\DefUnicode{σ}{\sigma}
\DefUnicode{ς}{\varsigma}
\DefUnicode{τ}{\tau}
\DefUnicode{υ}{\upsilon}
\DefUnicode{φ}{\varphi}
\DefUnicode{χ}{\chi}
\DefUnicode{ψ}{\psi}
\DefUnicode{ω}{\omega}

\DefUnicode{Γ}{\Gamma}
\DefUnicode{Δ}{\Delta}
\DefUnicode{Θ}{\Theta}
\DefUnicode{Λ}{\Lambda}
\DefUnicode{Ξ}{\Xi}
\DefUnicode{Π}{\Pi}
\DefUnicode{Σ}{\Sigma}
\DefUnicode{Φ}{\Phi}
\DefUnicode{Ψ}{\Psi}
\DefUnicode{Ω}{\Omega}


\DefUnicode{ℂ}{\CC}
\DefUnicode{ℍ}{\HH}
\DefUnicode{ℕ}{\NN}
\DefUnicode{ℙ}{\PP}
\DefUnicode{ℚ}{\QQ}
\DefUnicode{ℝ}{\RR}
\DefUnicode{ℤ}{\ZZ}

\DefUnicode{←}{\leftarrow}
\DefUnicode{→}{\rightarrow}
\DefUnicode{↛}{\nrightarrow}

\DefUnicode{∀}{\forall}
\DefUnicode{∃}{\exists}
\DefUnicode{∄}{\nexists}
\DefUnicode{∅}{\emptyset}
\DefUnicode{∈}{\in}
\DefUnicode{∉}{\notin}
\DefUnicode{∋}{\owns}
\DefUnicode{∌}{\notowns}
\DefUnicode{∏}{\prod}
\DefUnicode{∐}{\coprod}
\DefUnicode{∑}{\sum}
\DefUnicode{∘}{\circ}
\DefUnicode{∞}{\infty}
\DefUnicode{∧}{\wedge}
\DefUnicode{∨}{\vee}
\DefUnicode{∩}{\cap}
\DefUnicode{∪}{\cup}
\DefUnicode{≠}{\neq}
\DefUnicode{≤}{\leq}
\DefUnicode{≥}{\geq}
\DefUnicode{⊆}{\subseteq}
\DefUnicode{⊇}{\supseteq}
\DefUnicode{⊈}{\nsubseteq}
\DefUnicode{⊉}{\nsupseteq}
\DefUnicode{⊊}{\subsetneq}
\DefUnicode{⊋}{\supsetneq}
\DefUnicode{⋃}{\bigcup}
\DefUnicode{⋂}{\bigcap}

\DefUnicode{⋅}{\cdot}

\DefUnicode{⟨}{\langle}
\DefUnicode{⟩}{\rangle}

	\input{macros/local_labels}
	

\cdef \A {\mathcal{A}}
\cdef \B {\mathcal{B}}
\cdef \C {\mathcal{C}}
\cdef \D {\mathcal{D}}
\cdef \E {\mathcal{E}}
\cdef \F {\mathcal{F}}
\cdef \G {\mathcal{G}}
\cdef \H {\mathcal{H}}
\cdef \K {\mathcal{K}}
\cdef \M {\mathcal{M}}
\cdef \N {\mathcal{N}}
\cdef \O {\mathcal{O}}
\cdef \P {\mathcal{P}}
\cdef \R {\mathcal{R}}
\cdef \T {\mathcal{T}}
\cdef \U {\mathcal{U}}
\cdef \V {\mathcal{V}}
\cdef \W {\mathcal{W}}
\cdef \Y {\mathcal{Y}}

\cdef \NN {\mathbb{N}}
\cdef \RR {\mathbb{R}}
\cdef \CC {\mathbb{C}}


\cdef \set [1]{%
	\{#1\}%
}

\cdef \tuple [1]{%
	(#1)%
}


\cdef \card [1]{%
	\lvert #1\rvert
}

\cdef \powset [1]{%
	\mathcal{P}(#1)%
}

\cdef \subsets [2]{%
	[#2]^{#1}%
}

\cdef \disunion {%
	\sqcup
}

\cdef \DisUnion {%
	\bigsqcup
}

\cdef \continuum {%
	\mathfrak{c}%
}


\cdef \maps {%
	\colon
}

\cdef \into {%
	\hookrightarrow
}

\cdef \onto {%
	\twoheadrightarrow
}

\cdef \id {%
	\operatorname{id}%
}

\cdef \dom {%
	\operatorname{dom}%
}

\cdef \rng {%
	\operatorname{rng}%
}

\cdef \cod {%
	\operatorname{cod}%
}

\cdef \im [1]{%
	[#1]%
}

\cdef \inv {%
	^{-1}%
}

\cdef \preim [1]{%
	\inv\im{#1}%
}

\cdef \fiber [1]{%
	\inv(#1)%
}

\cdef \restr [1]{%
	\mathord{\upharpoonright}_{#1}%
}

\cdef \diag {%
	\mathbin\vartriangle
}

\cdef \Diag {%
	\bigtriangleup
}

\cdef \codiag {%
	\mathbin\triangledown
}

\cdef \CoDiag {%
	\bigtriangledown
}


\cdef \homeo {%
	\cong
}

\cdef \nothomeo {%
	\ncong
}

\cdef \clo [2][]{%
	\overline{#2}%
}

\cdef \cl {%
	\operatorname{cl}%
}

\cdef \topsum {%
	\oplus
}

\cdef \TopSum {%
	\sum
}


\cdef \dist {%
	\operatorname{d}%
}

\cdef \diam {%
	\operatorname{diam}%
}

\cdef \abs [1]{%
	\lvert #1\rvert
}


\cdef \Closed [1]{%
	\C l(#1)%
}

\cdef \Compacta [1]{%
	\K(#1)%
}

\cdef \Continua [2][]{%
	\C#1(#2)%
}

\cdef \ImageMap [1]{%
	#1^*%
}

\cdef \PreimageMap [1]{%
	#1^{-1**}%
}
\cdef \FiberMap [1]{%
	#1^{-1*}%
}


\cdef \half {%
	\frac{1}{2}%
}

\cdef \upset [1]{%
	#1^{\uparrow}%
}

\cdef \downset [1]{%
	#1^{\downarrow}%
}

\cdef \homeocopies [1]{%
	#1^{\homeo}%
}

\cdef \contimages [1]{%
	#1^{\twoheadrightarrow}%
}

\cdef \coll [1]{%
	[#1]%
}

\cdef \AllCompacta {%
	\mathbf{K}%
}

\cdef \AllContinua {%
	\mathbf{C}%
}

	
	\cdef \compactifiableHyper
		{\cite[Proposition~3.6]{compactifiable}}
	
	\cdef \compactifiableCompact
		{\cite[Theorem~3.13]{compactifiable}}
	
	\cdef \compactifiablePolish
		{\cite[Theorem~3.14]{compactifiable}}
	
	\cdef \compactifiableImplies
		{\cite[Corollary~3.17]{compactifiable}}
	
	\cdef \compactifiableExample
		{\cite[Example~4.15]{compactifiable}}
	
	\cdef \compactifiableInftyComponents
		{\cite[Example~4.18]{compactifiable}}
	
	\cdef \compactifiablePeano
		{\cite[Question~4.20]{compactifiable}}
	
	\cdef \compactifiableAnalytic
		{\cite[Theorem~4.26]{compactifiable}}

	\title{Borel complexity up to the equivalence}
	\author{Adam Bartoš \\
		\small
		Charles University, \\
		Faculty of Mathematics and Physics, \\
		Department of Mathematical Analysis, \\
		e-mail: \texttt{drekin@gmail.com}
	}
	\date{November 28, 2019}
	\maketitle
	
	\begin{abstract}
		We say that two classes of topological spaces are equivalent if each member of one class has a homeomorphic copy in the other class and vice versa. Usually when the Borel complexity of a class of metrizable compacta is considered, the class is realized as the subset of the hyperspace $\Compacta{[0, 1]^ω}$ containing all homeomorphic copies of members of the given class. We are rather interested in the lowest possible complexity among all equivalent realizations of the given class in the hyperspace.
		
		We recall that to every analytic subset of $\Compacta{[0, 1]^ω}$ there exists an equivalent $G_δ$ subset.
		Then we show that up to the equivalence open subsets of the hyperspace $\Compacta{[0, 1]^ω}$ correspond to countably many classes of metrizable compacta.
		Finally we use the structure of open subsets up to equivalence to prove that to every $F_σ$ subset of $\Compacta{[0, 1]^ω}$ there exists an equivalent closed subset.
		
		\begin{description}
			\item[Classification:] 
				54H05, 
				54B20, 
				54E45, 
				54F15. 
			
			\item[Keywords:] Borel hierarchy, complexity, homeomorphism equivalence, metrizable compactum, Polish space, hyperspace, Z-set, saturated family, compactifiable class, Polishable class.
		\end{description}
	\end{abstract}

	\linespread{1.2}\selectfont 

\section{Introduction}
	
	We denote that topological spaces $X$, $Y$ are \emph{homeomorphic} by $X \homeo Y$. This equivalence of topological spaces may be lifted to an equivalence of classes of topological spaces.
	We say that two classes $\C$ and $\D$ are \emph{equivalent} (and we also write $\C \homeo \D$) if every space in $\C$ is homeomorphic to a space in $\D$ and vice versa.
	This is the equivalence from the title.
	Given a class $\C$ we denote by $\homeocopies{\C}$ the class of all homeomorphic copies of members of $\C$.
	Clearly, this is the largest class equivalent to $\C$. We say that the class $\C$ is \emph{saturated} if $\C \homeo \homeocopies{\C}$.
	
	We denote the classes of all metrizable compacta and all metrizable continua by $\AllCompacta$ and $\AllContinua$, respectively.
	We are interested in the complexity of classes of metrizable compacta and continua, i.e. of subclasses of $\AllCompacta$ and $\AllContinua$.
	To express the complexity of a given class $\C$ using the Borel hierarchy, we first have to view the class as a subset of a Polish space.
	For this we use \emph{hyperspaces}.
	
	Let us recall the notation and basic properties of standard hyperspaces.
	For a topological space $X$ we denote the families of all compacta and continua (i.e. connected compacta) in $X$ including the empty set by $\Compacta{X}$ and $\Continua{X}$, respectively, and we endow the families with the \emph{Vietoris topology}.
	This is the topology generated by the sets of form $U^+ := \set{A: A ⊆ U}$ and $U^- := \set{A: A ∩ U ≠ ∅}$ for $U$ open in $X$.
	Clearly, $\Continua{X}$ is a subspace of $\Compacta{X}$.
	It is a closed subspace if $X$ is Hausdorff.
	
	If the space $X$ is metrizable with a metric $d$, the hyperspaces are metrizable with the induced \emph{Hausdorff metric} $d_H$. The distance $d_H(A, B)$ is the infimum of all values $ε > 0$ such that $A ⊆ N_ε(B)$ and $B ⊆ N_ε(A)$.
	Here, $N_ε(A)$ denotes the set $\set{x ∈ X: d(x, A) < ε} = ⋃_{x ∈ A} B(x, ε)$ where $B(x, ε)$ is the open ball of radius $ε$.
	To incorporate the empty set it makes sense to consider a bounded metric $d$ and to define $d_H(A, ∅)$ as the bound.
	
	Every continuous map $f\maps X \to Y$ between topological spaces induces the map $\ImageMap{f}\maps \Compacta{X} \to \Compacta{Y}$ defined by $\ImageMap{f}(A) := f\im{A}$. This map is also continuous. Moreover, if $f$ is an embedding or a homeomorphism, so is the map $\ImageMap{f}$.
	These properties are well known, and we summarize them in \compactifiableHyper.
	
	The \emph{Hilbert cube} $[0, 1]^ω$ is a universal space for all separable metrizable spaces, in particular, every metrizable compactum has a homeomorphic copy in $[0, 1]^ω$ regarded as an element of $\Compacta{[0, 1]^ω}$, which is itself a metrizable compactum. It is standard to view this space as the hyperspace of all metrizable compacta.
	For a class $\C ⊆ \AllCompacta$ we consider the collection of all families $\F ⊆ \Compacta{[0, 1]^ω}$ equivalent to $\C$, and we denote this collection by $\coll{\C}$.
	Note that this is the equivalence class of $\homeo$ restricted to $\powset{\Compacta{[0, 1]^ω}}$.
	Analogously to the saturated class we say that $\F ⊆ \Compacta{[0, 1]^ω}$ is a \emph{saturated family} if $\F = \homeocopies{\F} ∩ \Compacta{[0, 1]^ω}$.
	The collection $\coll{\C}$ has the largest element, namely the saturated family $\homeocopies{\C} ∩ \Compacta{[0, 1]^ω}$.
	Also, if $\F ∈ \coll{\C}$, then $\H ∈ \coll{\C}$ whenever $\F ⊆ \H ⊆ \max(\coll{\C})$.
	In particular, $\coll{\C}$ is stable under arbitrary unions.
	The minimal elements of $\coll{\C}$ are those families $\F ∈ \coll{\C}$ whose members are pairwise non-homeomorphic.
	
	Usually, when considering the complexity of a class $\C ⊆ \AllCompacta$, the class is identified with $\max(\coll{\C})$ and its complexity in $\Compacta{[0, 1]^ω}$ is considered.
	There are many results on complexity of $\max(\coll{\C})$, see for example the survey \cite{Marcone}.
	We are rather interested in the lowest complexity among families in $\coll{\C}$.
	This is rarely the complexity of the saturated family.
	For example, every singleton $\set{K} ⊆ \Compacta{[0, 1]^ω}$ is closed, but the corresponding saturated family is not unless $K$ is degenerate (see Observation~\ref{thm:singleton_saturation}).
	
	The reason we are interested in the lowest complexity among the members of $\coll{\C}$ for a class of metrizable compacta are the following notions introduced in \cite{compactifiable}.
	A class of topological spaces $\C$ is \emph{compactifiable} (or \emph{Polishable}) if there is a continuous map $q\maps A \to B$ between metrizable compacta (or Polish spaces) $A$, $B$ such that the family of fibers $\set{q\fiber{b}: b ∈ B}$ is equivalent to the given class $\C$.
	The map $q$ encodes how some representants of the members of $\C$ are disjointly composed together in one metrizable compactum (or Polish space) $A$. We call the resulting structure $\A(q\maps A \to B)$ a \emph{composition}.
	Clearly, if the class $\C$ is compactifiable (or Polishable), then it necessarily consists of metrizable compacta (or Polish spaces).
	
	We also define \emph{strongly compactifiable} and \emph{strongly Polishable classes} where the composition map $q$ additionally has to be closed and open.
	By \compactifiableImplies{} every compactifiable class is strongly Polishable. Therefore, we have the implications:
	\[
		\text{strongly compactifiable} \implies \text{compactifiable} \implies \text{strongly Polishable} \implies \text{Polishable}.
	\]
	The following theorems show that the strong notions are directly connected to hyperspaces, and since the definition of (strongly) compactifiable and Polishable classes is inherently up to the equivalence $\homeo$, we face the question of Borel complexity up to the equivalence.
	
	\begin{theorem}[\compactifiableCompact] \label{thm:strongly_compactifiable_characterization}
		The following conditions are equivalent for a class of topological spaces $\C$.
		\begin{enumerate}
			\item $\C$ is strongly compactifiable.
			\item There is a metrizable compactum $X$ and a closed family $\F ⊆ \Compacta{X}$ such that $\F \homeo \C$.
			\item There is a closed zero-dimensional disjoint family $\F ⊆ \Compacta{[0, 1]^ω}$ such that $\F \homeo \C$.
		\end{enumerate}
	\end{theorem}
	
	\begin{theorem}[\compactifiablePolish] \label{thm:strongly_Polishable_characterization}
		The following conditions are equivalent for a class of topological spaces $\C$.
		\begin{enumerate}
			\item $\C$ is a strongly Polishable class of compacta.
			\item There is a Polish space $X$ and an analytic family $\F ⊆ \Compacta{X}$ such that $\F \homeo \C$.
			\item There is a $G_δ$ zero-dimensional disjoint family $\F ⊆ \Compacta{[0, 1]^ω}$ such that $\F \homeo \C$.
			\item There is a closed zero-dimensional disjoint family $\F ⊆ \Compacta{(0, 1)^ω}$ such that $\F \homeo \C$.
		\end{enumerate}
	\end{theorem}
	
	So strong compactifiability correspond to existence of a closed equivalent subfamily of $\Compacta{[0, 1]^ω}$, and strong Polishability correspond to existence of an analytic or equivalently $G_δ$  equivalent subfamily of $\Compacta{[0, 1]^ω}$.
	The theorems are proved by translating back and forth between families in hyperspaces and compositions.
	As a byproduct, we obtain the following theorem. We include a sketch of a standalone proof that gathers all the translations needed together.
	
	\begin{theorem} \label{thm:analytic_G_delta}
		To every analytic family $\F ⊆ \Compacta{[0, 1]^ω}$ there exists an equivalent $G_δ$ family $\G ⊆ \Compacta{[0, 1]^ω}$.
		
		\begin{proof}
			Let $R := \set{\tuple{x, F} ∈ [0, 1]^ω × \F: x ∈ F}$ and let $π\maps R \to \F$ be the projection.
			Since the family $\F$ is analytic, there is a Polish space $B$ and a continuous surjection $f\maps B \to \F$.
			Let $A := \set{\tuple{x, b} ∈ [0, 1]^ω × B: x ∈ f(b)}$ and let $q\maps A \to B$ be the projection.
			The space $A$ is separable metrizable, and so there is an embedding $e\maps A \into [0, 1]^ω$.
			We put $\G := \set{e\im{q\fiber{b}}: b ∈ B}$.
			
			For every $b ∈ B$ we have $e\im{q\fiber{b}} \homeo q\fiber{b} = f(b) × \set{b} \homeo f(b) ∈ \F$, so $\G$ is equivalent to $\F$. The map $π$ is closed and open, and $q$ may be regarded as a pullback of $q$ along $f$. It follows that $q$ is also closed and open.
			We may suppose that $f$ also satisfies $\card{f\fiber{∅}} ≤ 1$.
			We obtain that $B$ is homeomorphic to $\set{q\fiber{b}: b ∈ B} ⊆ \Compacta{A}$, which is homeomorphic to $\G$ via $\ImageMap{e}$. Hence, $\G$ is Polish and so $G_δ$ in $\Compacta{[0, 1]^ω}$.
			For details see \cite{compactifiable}.
		\end{proof}
	\end{theorem}
	
	Let us note that for $σ$-ideals the previous theorem holds in a much stronger way.
	
	\begin{theorem}[{\cite[Theorem 11]{KLW_87}}] 
		Let $X$ be a metrizable compactum. Every analytic $σ$-ideal $\F ⊆ \Compacta{X}$ is in fact $G_δ$.
	\end{theorem}
	
	In this paper we analyze the remaining complexities, namely clopen, open, and $F_σ$ subsets of $\Compacta{[0, 1]^ω}$.
	The situation with clopen subsets is quite simple. It is well-known that the hyperspaces $\Compacta{X} \setminus \set{∅}$ and $\Continua{X} \setminus \set{∅}$ are connected for any connected space $X$ (see for example \cite[Exercises 4.32 and 5.25]{Nadler}). Hence, we obtain the following proposition.
	
	\begin{proposition} \label{thm:clopen_classes}
		There are exactly four clopen subsets of $\Compacta{[0, 1]^ω}$: $∅$, $\set{∅}$, $\Compacta{[0, 1]^ω} \setminus \set{∅}$, $\Compacta{[0, 1]^ω}$.
		Hence, there are only four corresponding classes: $∅$, $\set{∅}$, $\AllCompacta \setminus \set{∅}$, $\AllCompacta$.
		Similarly, there are exactly four clopen subsets of $\Continua{[0, 1]^ω}$: $∅$, $\set{∅}$, $\Continua{[0, 1]^ω} \setminus \set{∅}$, $\Continua{[0, 1]^ω}$, and four corresponding classes of continua: $∅$, $\set{∅}$, $\AllContinua \setminus \set{∅}$, $\AllContinua$.
	\end{proposition}
	
	The situation with open and $F_σ$ families is more involved and is the subject of the next sections.
	In the second section we prove that every open subset of $\Compacta{[0, 1]^ω}$ is equivalent to one of countably many saturated open subfamilies of the hyperspace (Theorem~\ref{thm:open_classes}).
	In the third section we show that every $F_σ$ subset of $\Compacta{[0, 1]^ω}$ is equivalent to a closed subset (Theorem~\ref{thm:F_sigma_classes}).
	In the fourth section we gather some observations on saturated and so-called type-saturated classes and families.

\section{Open classes}
	
	Now let us look at open subsets of $\Compacta{[0, 1]^ω}$ up to the equivalence.
	First, we shall consider the following rough classification of metrizable compacta.
	
	\begin{definition}
		Let $X$ be a metrizable compactum.
		\begin{itemize}
			\item By $m(X)$ we denote the number of all connected components.
			By $n(X)$ we denote the number of all nondegenerate connected components.
			
			\item Let $T$ denote the set of all \emph{finite types} $\set{\tuple{m, n}: m ≥ n ∈ ω}$, and let $T_+$ denote the set of all \emph{positive finite types} $\set{\tuple{m, n} ∈ T: m > 0}$.
			
			\item We define the \emph{type function} $t\maps \AllCompacta \to T ∪ \set{∞}$ by $t(X) := \tuple{m(X), n(X)}$ if $m(X) < ω$, $∞$ otherwise.
				Clearly, the type function is onto.
			
			\item We define a partial order $≤$ on $T ∪ \set{∞}$: $\tuple{0, 0}$ is not comparable with anything; 
				$T_+$ is ordered by the product order, i.e. $\tuple{m_1, n_1} ≤ \tuple{m_2, n_2}$ if and only if $m_1 ≤ m_2$ and $n_1 ≤ n_2$;
				and $∞ ≥ t$ for every $t ∈ T_+$.
			
			\item We define the \emph{principal upper class} $\U_t := \set{X ∈ \AllCompacta: t(X) ≥ t}$ for every $t ∈ T ∪ \set{∞}$.
				Since the type function is onto, we have $t = \min\set{t(X): X ∈ \U_t}$ for every $t ∈ T ∪ \set{∞}$, and so $t_1 ≤ t_2 \iff \U_{t_1} ⊇ \U_{t_2}$ for every $t_1, t_2 ∈ T ∪ \set{∞}$.
		\end{itemize}
	\end{definition}
	
	\begin{example}
		We have the following examples of principal upper classes.
		\begin{itemize}
			\item $\U_{m, 0}$ is the class of all metrizable compacta with at least $m$ components.
			\item $\U_{m, 0} ∪ \U_{1, 1}$ is the class of all metrizable compacta with at least $m$ points.
			\item $\U_{2, 0} ∪ \U_{1, 1}$ is the class of all nondegenerate metrizable compacta.
			\item $\U_{1, 1}$ is the class of all infinite metrizable compacta.
			\item $\U_{1, 0}$ is the class of all nonempty metrizable compacta, i.e.\ $\AllCompacta \setminus \set{∅}$.
			\item $\U_{0, 0} = \set{∅}$ and $\U_{0, 0} ∪ \U_{1, 0} = \AllCompacta$.
		\end{itemize}
	\end{example}
	
	We will show that open subsets of $\Compacta{[0, 1]^ω}$ are equivalent to some unions of principal upper classes.
	Since the finite spaces are dense in $\Compacta{[0, 1]^ω}$, not every principal upper class is open.
	However, this is essentially the only obstacle.
	That is why we define \emph{nice} sets of types.
	
	\begin{definition}
		Let $R ⊆ T ∪ \set{∞}$.
		\begin{itemize}
			\item We say that $R$ is \emph{nice} if $\tuple{m, 0} ∈ R$ for some $m > 0$ whenever $R ∩ (T_+ ∪ \set{∞}) ≠ ∅$.
				This holds if and only if $⋃_{t ∈ R} \U_t$ contains a nonempty finite space whenever it contains a nonempty space.
			
			\item We say that $R$ is an \emph{antichain} if it is pairwise $≤$-incomparable.
				Note that every antichain is finite, and that no nice antichain contains $∞$.
			
			\item By $A(R)$ we denote the set of all $≤$-minimal elements of $R$.
				Note that this is the only antichain $A$ such that $⋃_{t ∈ A} \U_t = ⋃_{t ∈ R} \U_t$.
				It follows that $A(R)$ is nice if and only if $R$ is nice.
		\end{itemize}
	\end{definition}
	
	Eventually, we will show that open subsets of $\Compacta{[0, 1]^ω}$ correspond to nice antichains in $T$ (Theorem~\ref{thm:open_classes}), but first we determine which unions of principal upper classes are open.
	
	\begin{definition}
		For every finite function $s\maps I \to ℕ_+$, where $ℕ_+$ denotes the set of all positive integers, we define the \emph{special open class} $\O_s$ of all metrizable compacta $K$ having a clopen decomposition $\set{K_i: i ∈ I}$ such that $\card{K_i} ≥ s(i)$ for every $i ∈ I$.
		
		Moreover, let $X$ be a metrizable space, and let $\U ⊆ \Compacta{X}$ be open.
		We say that $\U$ is \emph{of the shape $s$} if there are disjoint open sets $U_i ⊆ X$, $i ∈ I$, and for every $i ∈ I$ there are disjoint open sets $V_{i, j} ⊆ U_i$, $j < s(i)$, such that $\U = (⋃_{i ∈ I} U_i)^+ ∩ ⋂_{i ∈ I, j < s(i)} V_{i, j}^-$.
		We say that $\U$ is \emph{exactly of the shape $s$} if moreover the set $U_i^+ ∩ ⋂_{j < s(i)} V_{i, j}^-$ contains a connected space for every $i ∈ I$.
		
		By $n(s)$ we denote $\card{\set{i ∈ I: s(i) > 1}}$.
		To every type $t ∈ T ∪ \set{∞}$ we associate a set of finite functions $S_t$.
		If $t = \tuple{m, n}$, we put $S_t := \set{s\maps m \to ℕ_+: n(s) ≤ n}$,
		if $t = ∞$, we put $S_t := \set{s\maps m \to ℕ_+: m > 0}$.
	\end{definition}
	
	\begin{observation} \label{thm:special_open_class}
		Let $s\maps I \to ℕ_+$ be a finite function, let $X$ be a metrizable space, and let $K ∈ \Compacta{X}$.
		$K$ has a neighborhood of the shape $s$ in $\Compacta{X}$ if and only if $K ∈ \O_s$.
		It follows that $\O_s ∩ \Compacta{X}$ is open.
	\end{observation}
	
	\begin{observation} \label{thm:quasi-components}
		Let $s\maps I \to ℕ_+$ be a nonempty finite function and let be $K$ a metrizable compactum.
		If there are pairwise disjoint sets $A_i ⊆ K$, $i ∈ I$, such that for every $i ∈ I$ either $A_i$ is a nondegenerate component of $K$ or $A_i$ is the union of $s(i)$-many components, then $K ∈ \O_s$.
		
		This is because the components and the quasi-components are the same and we have used only finitely many components when building the sets $A_i$, and hence there is a clopen decomposition $\set{K_i: i ∈ I}$ of $K$ such that $A_i ⊆ K_i$ for every $i ∈ I$.
		Also, every nondegenerate component is infinite, so $\card{K_i} ≥ \card{A_i} ≥ s(i)$ for every $i ∈ I$.
	\end{observation}
	
	Note that each antichain in $T_+$ is of the form $\set{\tuple{m + ∑_{i < j} Δm_i,\ n - ∑_{i < j} Δn_i}: j ≤ k}$ for some $\set{Δm_i, Δn_i: i < k} ⊆ ℕ_+$, and it is nice if and only if $∑_{i < k} Δn_i = n$, so the last member is $\tuple{m + ∑_{i < k} Δm_i, 0}$.
	The next proposition says that each special open class $\O_s$ corresponds to such nice antichain additionally satisfying that each $Δn_i$ is $1$ and that the sequence $\tuple{Δm_i: i < k}$ is increasing.
	
	\begin{proposition} \label{thm:special_antichain}
		Let $s\maps I \to ℕ_+$ be a finite function. We have $\O_s = ⋃_{t ∈ R_s} \U_t$ where $R_s$ is a nice antichain in $T$ defined as follows.
		
		Let $\tuple{i_k: k < \card{I}}$ be an enumeration of $I$ such that the map $k \mapsto s(i_k)$ is increasing.
		For every $n ≤ n(s)$ let us consider the type $t_{s, n} := \tuple{n + ∑_{k < \card{I} - n} s(i_k), n}$.
		In particular, $t_{s, 0} = \tuple{∑_{i ∈ I} s(i), 0}$ and $t_{s, n(s)} = \tuple{\card{I}, n(s)}$.
		We put $R_s := \set{t_{s, n}: n ≤ n(s)}$.
		
	
		\begin{proof}
			First, if $s = ∅$, we have $\O_s = \set{∅} = \U_{0, 0} = \U_{t_{s, 0}}$, so we may suppose that $s ≠ ∅$.
			
			If $K ∈ \O_s$, then it has a clopen decomposition $\set{K_i: i ∈ I}$ such that for every $i ∈ I$ we have $\card{K_i} ≥ s(i)$.
			Let $J := \set{i ∈ I: s(i) > 1$ and $K_i$ contains a nondegenerate component$}$ and $n := \card{J}$.
			Clearly, $n ≤ n(s)$.
			We have $∑_{k < \card{I} - n} s(i_k) ≤ ∑_{i ∈ I \setminus J} s(i)$ since the map $k \mapsto s(i_k)$ is increasing.
			Therefore, $t(K) ≥ \tuple{\card{J} + ∑_{i ∈ I \setminus J} s(i), \card{J}} ≥ \tuple{n + ∑_{k < \card{I} - n} s(i_k), n} = t_{s, n}$.
			It follows that $\O_s ⊆ ⋃_{t ∈ R_s} \U_t$.
			
			On the other hand, if $K ∈ \U_{t_{s, n}}$ for some $n ≤ n(s)$, then $K$ has at least $n + ∑_{k < \card{I} - n} s(i_k)$ components at least $n$ of which are nondegenerate.
			Hence, we may find disjoint sets $A_i ⊆ K$, $i ∈ I$, such that $A_{i_k}$ is a nondegenerate component if $k ≥ \card{I} - n$ and $A_{i_k}$ is the union of $s(i_k)$ components if $k < \card{I} - n$.
			From Observation~\ref{thm:quasi-components} it follows that $K ∈ \O_s$, and so $⋃_{t ∈ R_s} \U_t ⊆ \O_s$.
		\end{proof}
	\end{proposition}
	
	\begin{example}
		We have the following examples of special open classes.
		\begin{itemize}
			\item $\O_∅ = \U_{0, 0} = \set{∅}$ is the empty space class.
			\item $\O_{\tuple{1}} = \U_{1, 0} = \AllCompacta \setminus \set{∅}$ is the class of all nonempty metrizable compacta.
			\item $\O_{\tuple{2}} = \U_{1, 1} ∪ \U_{2, 0}$ is the class of all nondegenerate metrizable compacta.
			\item $\O_{\tuple{m}} = \U_{1, 1} ∪ \U_{m, 0}$ is the class of all metrizable compacta with at least $m$ points.
			\item $\O_{\tuple{1: i < m}} = \U_{m, 0}$ is the class of all metrizable compacta with at least $m$ components.
			\item $\O_{\tuple{1, 1, 1, 2, 3, 4}} = \U_{6, 3} ∪ \U_{7, 2} ∪ \U_{9, 1} ∪ \U_{12, 0}$.
		\end{itemize}
	\end{example}
	
	\begin{corollary} \label{thm:open_sandwich}
		For every $t ∈ T ∪ \set{∞}$ and every $m ∈ ℕ_+$ there is $s_{t, m} ∈ S_t$ such that $\U_t ⊆ \O_{s_{t, m}} ⊆ \U_t ∪ \U_{m, 0}$.
		
		\begin{proof}
			For $t = ∞$ we simply put $s_{t, m} := \tuple{1: i < m}$ so $\O_{s_{t, m}} = \U_{m, 0}$.
			For $t = \tuple{m', n'} ∈ T$ we define $s_{t, m} = s$ as a function with domain $m'$ taking the value $m$ $n'$ times and the value $1$ $m' - n'$ times.
			By Proposition~\ref{thm:special_antichain} we have $\O_{s_{t, m}} = ⋃_{n ≤ n'} \U_{t_{s, n}}$ and $t_{s, n} = \tuple{n + (m' - n') + (n' - n) ⋅ m,\ n}$.
			Hence, for $n = n'$ we obtain $\U_{t_{s, n}} = \U_t$ and for $n' - n > 0$ the first item is at least $m$, so $\U_{t_{s, n}} ⊆ \U_{m, 0}$.
		\end{proof}
	\end{corollary}
	
	\begin{proposition} \label{thm:principal_upper_class}
		For every $t ∈ T ∪ \set{∞}$ we have $\U_t = ⋂_{s ∈ S_t} \O_s$.
		In particular, $\U_t ∩ \Compacta{X}$ is $G_δ$ for every metrizable space $X$, so every principal upper class is strongly Polishable.
		It also follows that $\U_{t'} ⊆ \O_s$ for every $t' ≥ t$ and $s ∈ S_t$.
		
		\begin{proof}
			First let us show that $\U_t ⊆ ⋂_{s ∈ S_t} \O_s$, so let $K ∈ \U_t$ and $s ∈ S_t$.
			If $t = \tuple{m, n} ∈ T$, then $K$ has a clopen decomposition $\set{K_i: i < m}$ into components. Since $n(s) ≤ n$, we may choose the enumeration such that $K_i$ is nondegenerate whenever $s(i) > 1$. Since nondegenerate components are infinite, we have $\card{K_i} ≥ s(i)$ for every $i < m$.
			If $t = ∞$, then $K$ has infinitely many components, so we may find suitable sets $A_i$ and use Observation~\ref{thm:quasi-components}.
			In both cases we have $K ∈ \O_s$.
			
			Now, $\U_t ⊇ ⋂_{s ∈ S_t} \O_s$.
			If $t ≤ ∞$, then for every $m > 0$ we take $s_{t, m} ∈ S_t$ from Corollary~\ref{thm:open_sandwich}, and we have 
			$\U_t ⊆ ⋂_{m ∈ ℕ_+} \O_{s_{t, m}} ⊆ \U_t ∪ ⋂_{m ∈ ℕ_+} \U_{m, 0} = \U_t ∪ \U_∞ = \U_t$.
			Otherwise, $t = \tuple{0, 0}$ and $\U_t = \set{∅} = \O_∅$.
		\end{proof}
	\end{proposition}
	
	\begin{proposition} \label{thm:open_unions}
		Let $R ⊆ T ∪ \set{∞}$. The set $⋃_{t ∈ R} \U_t ∩ \Compacta{[0, 1]^ω}$ is open if and only if $R$ is nice.
		
		\begin{proof}
			First, suppose that $R$ is nice.
			Let $t ∈ R$.
			If $t = \tuple{0, 0}$ we put $s_t := ∅$ and we have $\U_t = \O_{s_t}$.
			Otherwise, there is $m > 0$ such that $\tuple{m, 0} ∈ R$, and we put $s_t := s_{t, m}$ from Corollary~\ref{thm:open_sandwich}, so $\U_t ⊆ \O_{s_t} ⊆ \U_t ∪ \U_{m, 0}$.
			Altogether, we have $⋃_{t ∈ R} \U_t = ⋃_{t ∈ R} \O_{s_t}$, which has open intersection with $\Compacta{[0, 1]^ω}$ by Observation~\ref{thm:special_open_class}.
			
			On the other hand, if $\U := ⋃_{t ∈ R} \U_t ∩ \Compacta{[0, 1]^ω}$ is open and $R$ meets $T_+ ∪ \set{∞}$, we have $\U \setminus \set{∅} ≠ ∅$.
			Since finite sets are dense, there is a finite set $F ∈ \U \setminus \set{∅}$, and there is some $t ∈ R$ such that $F ∈ \U_t$. Since $F$ is finite and nonempty, we have $t = \tuple{m, 0}$ for some $m > 0$, so $R$ is nice.
		\end{proof}
	\end{proposition}
	
	The previous propositions regarding the properties of principal upper classes and special open classes would hold as well in the realm of Hausdorff compacta instead of metrizable compacta. Hausdorffness is needed so that components and quasi-components are the same in compacta and that nondegenerate connected spaces are infinite.
	
	We have shown that open unions of principal upper classes are exactly unions over nice antichains. Now we show that every open subset of $\Compacta{[0, 1]^ω}$ is equivalent to such union.
	
	The following lemma and also the more general fact that the set of all homeomorphic copies of any nondegenerate continuum is dense in $\Continua{[0, 1]^ω}$ should be known, but for the reader’s convenience a short proof is included.
	
	\begin{lemma} \label{thm:Hilbert_cube_dense}
		The set of all homeomorphic copies of $[0, 1]^ω$ is dense in $\Continua{[0, 1]^ω} \setminus \set{∅}$.
		
		\begin{proof}
			Let $U^+ ∩ ⋂_{i < n} V_i^-$ be a basic neighborhood of a nonempty continuum $C ⊆ [0, 1]^ω$.
			Since $C$ is connected and $[0, 1]^ω$ is locally path-connected, we may suppose that the set $U$ is path-connected.
			For $i < n$ let $y_i ∈ U ∩ V_i$, and let $Y$ be the union of finitely many paths in $U$ connecting the points $y_i$.
			There is some $ε > 0$ such that $N_ε(Y) ⊆ U$.
			Let $f\maps [0, 1]^ω \to Y$ be a continuous surjection, for every $i < n$ let $x_i ∈ [0, 1]^ω$ be such that $f(x_i) = y_i$, and let $A := \set{x_i: i < n}$.
			By the Mapping Replacement Theorem \cite[5.3.11]{Mill} there is a Z-embedding $g\maps [0, 1]^ω \to [0, 1]^ω$ such that $g\restr{A} = f\restr{A}$ and $d(g, f) < ε$. Therefore, $[0, 1]^ω \homeo \rng(g) ∈ U^+ ∩ ⋂_{i < n} V_i^-$.
		\end{proof}
	\end{lemma}
	
	\begin{lemma} \label{thm:big_in_Hilbert_cube}
		Let $F ⊆ [0, 1]^ω$ be a finite set. For every separable metrizable space $X$ such that $\card{X} ≥ \card{F}$ there exists an embedding $f\maps X \into [0, 1]^ω$ such that $F ⊆ f\im{X}$.
		
		\begin{proof}
			Since $X$ is separable metrizable, we may suppose that $X ⊆ [0, 1]^ω$. Since $\card{X} ≥ \card{F}$, there is a bijection $h\maps H \to F$ for some $H ⊆ X$. The map $h$ is a homeomorphism of Z-sets in $[0, 1]^ω$, so by \cite[Theorem~5.3.7]{Mill} it can be extended to a homeomorphism $\bar{h}\maps [0, 1]^ω \to [0, 1]^ω$. The restriction $\bar{h}\restr{X}$ is the desired embedding.
		\end{proof}
	\end{lemma}
	
	\begin{proposition} \label{thm:exact_shape_contains_space}
		Let $s\maps I \to ℕ_+$ be a finite function.
		For every compactum $X ∈ \O_s$ and every open set $\U ⊆ \Compacta{[0, 1]^ω}$ exactly of the shape $s$ there is a compactum $Y ∈ \U$ homeomorphic to $X$.
		
		\begin{proof}
			Let $\set{U_i, V_{i, j}: i ∈ I, j < s(i)}$ be the open subsets of $[0, 1]^ω$ witnessing that $\U$ is exactly of the shape $s$, and let $\set{X_i: i ∈ I}$ be a clopen decomposition of $X$ such that $\card{X_i} ≥ s(i)$ for every $i ∈ I$.
			Let $i ∈ I$.
			Since $U_i^+ ∩ ⋂_{j < s(i)} V_{i, j}^-$ contains a connected space, it also contains a space $Q_i \homeo [0, 1]^ω$ by Lemma~\ref{thm:Hilbert_cube_dense}.
			Let $F_i ⊆ Q_i$ be such that $\card{F_i} = s(i)$ and $F_i ∩ V_{i, j} ≠ ∅$ for every $j < s(i)$.
			By Lemma~\ref{thm:big_in_Hilbert_cube} there is a copy $Y_i \homeo X_i$ such that $F_i ⊆ Y_i ⊆ Q_i$.
			Hence, $Y_i ∈ U_i^+ ∩ ⋂_{j < s(i)} V_{i, j}^-$.
			Altogether we have $X \homeo Y := ⋃_{i ∈ I} Y_i ∈ \U$.
		\end{proof}
	\end{proposition}
	
	\begin{lemma} \label{thm:exact_shape_basis}
		Let $t ∈ T ∪ \set{∞}$.
		Every $K ∈ \U_t ∩ \Compacta{X}$ for any metrizable space $X$ has a neighborhood basis such that for every basic set $\U$ there is $s ∈ S_t$ such that $\U$ is exactly of the shape $s$.
		
		\begin{proof}
			Let $\V ⊆ \Compacta{X}$ be any neighborhood of $K$. 
			Without loss of generality $\V$ is of the form $V^+ ∩ ⋂\set{W^-: W ∈ \W}$ for some open set $V ⊆ X$ and a finite family of open sets $\W$.
			
			If $t = ∞$, let $\set{C_i: i ∈ I}$ be a finite collection of distinct components of $K$ such that every $W ∈ \W$ meets some of them, and let $\set{K_i: i ∈ I}$ be a clopen decomposition of $K$ such that $C_i ⊆ K_i$ for every $i ∈ I$. Such sets $K_i$ exist since components of $K$ are the quasi-components.
			If $t = \tuple{m, n} ∈ T$, let $\set{C_i = K_i: i ∈ I = m}$ be the enumeration of all components of $K$.
			
			For every $i ∈ I$ let $F_i := \set{x_{i, j}: j < s(i)} ⊆ C_i$ be a nonempty finite set of minimal size such that $F_i ∩ W ≠ ∅$ for every $W ∈ \W ∩ C_i^-$.
			This defines a function $s\maps I \to ℕ_+$.
			For every $i ∈ I$ we have $s(i) ≤ \card{C_i}$, and so $n(s) ≤ n$ if $t = \tuple{m, n}$. Hence, $s ∈ S_t$.
			
			Since the set $I$ is finite, there are disjoint open sets $U_i ⊆ V$, $i ∈ I$, such that $K_i ⊆ U_i$,
			and for every $i ∈ I$ there are disjoint open sets $U_{i, j} ⊆ U_i$, $j < s(i)$, such that $x_{i, j} ∈ U_{i, j} ⊆ ⋂\set{W ∈ \W: x_{i, j} ∈ W}$.
			We put $\U := (⋃_{i ∈ I} U_i)^+ ∩ ⋂_{i ∈ I, j < s(i)} U_{i, j}^-$ and $\U_i := U_i^+ ∩ ⋂_{j < s(i)} U_{i, j}^-$ for every $i ∈ I$.
			Since $⋃_{i ∈ I} U_i ⊆ V$ and for every $W ∈ \W$ there is $i ∈ I$ and $j < s(i)$ such that $U_{i, j} ⊆ W$, we have $\U ⊆ \V$.
			Since $C_i, K_i ∈ \U_i$ for every $i ∈ I$, we have that $\U$ is exactly of the shape $s$ and $K ∈ \U$.
		\end{proof}
	\end{lemma}
	
	\begin{proposition} \label{thm:neighborhood_types}
		Let $X, Y ∈ \Compacta{[0, 1]^ω}$.
		A homeomorphic copy of $Y$ is contained in every neighborhood of $X$ if and only if $t(Y) ≥ t(X)$.
		
		\begin{proof}
			\BackwardImplication: Suppose that $t(Y) ≥ t(X)$ and let $\U$ be a neighborhood of $X$.
			By Lemma~\ref{thm:exact_shape_basis} we may suppose that $\U$ is exactly of the shape $s$ for some $s ∈ S_{t(X)}$.
			By Proposition~\ref{thm:principal_upper_class} we have $Y ∈ \U_{t(Y)} ⊆ \U_{t(X)} ⊆ \O_s$.
			Finally, by Proposition~\ref{thm:exact_shape_contains_space}, there is a space $Y' ∈ \U$ homeomorphic to $Y$.
			
			\ForwardImplication: Suppose that $t(Y) \ngeq t(X)$.
			We have $Y ∉ \U_{t(X)} = ⋂_{s ∈ S_{t(X)}} \O_s$ by Proposition~\ref{thm:principal_upper_class}.
			Hence, there is some $s ∈ S_{t(X)}$ such that $Y ∉ \O_s ∩ \Compacta{[0, 1]^ω} ∋ X$.
			Since $\O_s$ is closed under homeomorphic copies, we are done.
		\end{proof}
	\end{proposition}
	
	\begin{definition} \label{def:open_class}
		By $\R$ we denote the countable set of all nice antichains of $T ∪ \set{∞}$.
		For every $R ∈ \R$ we define the \emph{open class} $\O_R := ⋃_{t ∈ R} \U_t$.
		Proposition~\ref{thm:special_antichain} says that every special open class is an open class, namely $\O_s = \O_{R_s}$ for every finite $s\maps I \to ℕ_+$.
	\end{definition}
	
	\begin{theorem} \label{thm:open_classes}
		For every open $\U ⊆ \Compacta{[0, 1]^ω}$ there exists exactly one $R ∈ \R$ such that $\U \homeo \O_R$.
		On the other hand, for every $R ∈ \R$ we have $\O_R \homeo \O_R ∩ \Compacta{[0, 1]^ω}$, which is open.
		
		\begin{proof}
			By Proposition~\ref{thm:neighborhood_types} and by universality of $\Compacta{[0, 1]^ω}$ we have $\U \homeo ⋃_{X ∈ \U} \U_{t(X)}$. We put $R := A(\set{t(X): X ∈ \U})$. Since $\U$ is open, it contains a nonempty finite space whenever it contains a nonempty space. Therefore, $R$ is nice and $\U \homeo \O_R$.
			
			Clearly, if $R ≠ R' ∈ \R$, there is a type $t ∈ T$ that is above some member of $R$ and above no member of $R'$ or the other way around.
			Any metrizable compactum $X$ of type $t$ satisfies $X ∈ (\O_R \setminus \O_{R'}) ∪ (\O_{R'} \setminus \O_R)$, and hence $\O_R \nothomeo \O_{R'}$.
			
			On the other hand, let $R ∈ \R$.
			$\O_R ∩ \Compacta{[0, 1]^ω}$ is open by Proposition~\ref{thm:open_unions}, and $\O_R \homeo \O_R ∩ \Compacta{[0, 1]^ω}$ since $\Compacta{[0, 1]^ω}$ is universal for metrizable compacta.
		\end{proof}
	\end{theorem}
	
	\begin{corollary}
		There are exactly six nonequivalent classes corresponding to open subsets of $\Continua{[0, 1]^ω}$.
		Besides the four clopen classes $∅$, $\set{∅}$, $\AllContinua \setminus \set{∅}$, and $\AllContinua$, there is the class of all nondegenerate continua $\U_{1, 1} ∩ \AllContinua$ and the class $(\U_{1, 1} ∪ \U_{0, 0}) ∩ \AllContinua = (\U_{1, 1} ∩ \AllContinua) ∪ \set{∅}$.
		
		\begin{proof}
			Every open subset $\V$ of $\Continua{[0, 1]^ω}$ is of form $\U ∩ \AllContinua$ where $\U$ is open in $\Compacta{[0, 1]^ω}$.
			By Theorem~\ref{thm:open_classes} we have $\U \homeo \O_R$ for some nice antichain $R$, and hence $\V \homeo \O_R ∩ \AllContinua$.
			Since $\U_{2, 0} ∩ \AllContinua = ∅$, open subsets of $\Continua{[0, 1]^ω}$ are equivalent to classes $⋃_{t ∈ R} \U_t ∩ \AllContinua$ where $R$ is any antichain in $\set{\tuple{0, 0}, \tuple{1, 0}, \tuple{1, 1}}$.
			These are the six declared classes.
		\end{proof}
	\end{corollary}

\section{Countable unions of strongly compactifiable classes}
	
	In this section we show that every $F_σ$ subset of $\Compacta{[0, 1]^ω}$ is equivalent to a closed subset, or equivalently, that strongly compactifiable classes are stable under countable unions.
	But first we have to improve several results from the previous section.
	
	\begin{lemma} \label{thm:closed_intersecting}
		Let $X$ be a metrizable space and let $\F ⊆ \Compacta{X}$ be a compact family.
		For every open set $U ⊆ X$ such that $\F ⊆ U^-$ there exists a closed set $A ⊆ X$ such that $A ⊆ U$ and $\F ⊆ A^-$.
		
		\begin{proof}
			Let $d$ be a compatible metric on $X$.
			For every $F ∈ \F$ there is $x_F ∈ F$ and $δ_F > 0$ such that $B(x_F, δ_F) ⊆ U$.
			Since $\F$ is compact, there is a finite collection $\H ⊆ \F$ such that $\F ⊆ ⋃_{H ∈ \H} B(x_H, δ_H / 2)^-$.
			Hence, for every $F ∈ \F$ there is $H_F ∈ \H$ and $y_F ∈ F ∩ B(x_{H_F}, δ_{H_F} / 2)$.
			We put $Y := \set{y_F: F ∈ \F}$ and $δ := \min\set{δ_H / 2: H ∈ \H}$.
			For every $F ∈ \F$ we have that $B(y_F, δ) ⊆ B(x_{H_F}, δ_{H_F}) ⊆ U$.
			Therefore, $d(Y, X \setminus U) ≥ δ$ and $A := \clo{Y} ⊆ U$.
		\end{proof}
	\end{lemma}
	
	\begin{lemma} \label{thm:bigger_in_Hilbert_cube}
		Let $X$ be a separable metrizable space, let $J$ be finite, and let $F_j ⊆ X$, $j ∈ J$, be disjoint compact sets.
		Let $V_j ⊆ [0, 1]^ω$, $j ∈ J$, be disjoint nonempty open sets.
		There is an embedding $f\maps X \into [0, 1]^ω$ such that $f\im{F_j} ⊆ V_j$ for every $j ∈ J$.
		
		\begin{proof}
			There exists a Z-set $Q ∈ ⋂_{j ∈ J} V_j^-$ such that $Q \homeo [0, 1]^ω$. This follows from \cite[Lemma~5.1.3]{Mill} since there is $n ∈ ω$ such that every set $V_j$ contains a point $x_j$ such that $π_n(x_j) = 1$.
			Also, by Lemma~\ref{thm:Hilbert_cube_dense} there are sets $Q_j ⊆ Q ∩ V_j$, $j ∈ J$, such that $Q_j \homeo [0, 1]^ω$ for every $j ∈ J$.
			
			Since $X$ is separable metrizable, we may suppose that $X ⊆ Q$.
			There are homeomorphisms $h_j\maps F_j \to H_j ⊆ Q_j$ for $j ∈ J$.
			The map $h := ⋃_{j ∈ J} h_j$ is a homeomorphism of Z-sets in $[0, 1]^ω$ since $⋃_{j ∈ J} F_j$ and $⋃_{j ∈ J} H_j$ are closed subsets of the Z-set $Q$.
			By \cite[Theorem~5.3.7]{Mill} the map $h$ can be extended to a homeomorphism $\bar{h}\maps [0, 1]^ω \to [0, 1]^ω$.
			The restriction $\bar{h}\restr{X}$ is the desired embedding.
		\end{proof}
	\end{lemma}
	
	\begin{proposition} \label{thm:exact_shape_transport}
		Let $s\maps I \to ℕ_+$ be a finite function, let $\U ⊆ \Compacta{[0, 1]^ω}$ be an open set exactly of the shape $s$, and let $\V ⊆ \Compacta{X}$ be an open set of the shape $s$ for some metrizable space $X$.
		For every compact family $\H ⊆ \V$ there is a compact family $\F ⊆ \U$ and a homeomorphism $Φ\maps \H \to \F$ such that $Φ(H) \homeo H$ for every $H ∈ \H$.
		
		\begin{proof}
			Let $\set{U_i, U_{i, j}: i ∈ I, j < s(i)}$ be open subsets of $[0, 1]^ω$ witnessing that $\U$ is exactly of the shape $s$,
			let $\set{V_i, V_{i, j}: i ∈ I, j < s(i)}$ be open subsets of $X$ witnessing that $\V$ is of the shape $s$,
			and let $\H ⊆ \V$ be a compact family.
			We fix $i ∈ I$ and put $\U_i := U_i^+ ∩ ⋂_{j < s(i)} U_{i, j}^-$.
			Since $\U_i$ contains a connected space, it also contains a space $Q_i \homeo [0, 1]^ω$ by Lemma~\ref{thm:Hilbert_cube_dense}.
			For every $j < s(i)$ there is a compact set $A_{i, j} ⊆ V_{i, j}$ such that $\H ⊆ A_{i, j}^-$ (Lemma~\ref{thm:closed_intersecting}).
			By Lemma~\ref{thm:bigger_in_Hilbert_cube} there is an embedding $e_i\maps [0, 1]^ω \to Q_i$ such that $e_i\im{A_{i, j}} ⊆ U_{i, j}$ for every $j < s(i)$.
			
			For every $i ∈ I$ we have the homeomorphism $h_i := e_i\restr{V_i}\maps V_i \to \rng(e_i) ⊆ Q_i$. Since the families $\set{V_i: i ∈ I}$ and $\set{\rng(e_i): i ∈ I}$ are separated, the map $h := ⋃_{i ∈ I} h_i\maps ⋃_{i ∈ I} V_i \to ⋃_{i ∈ I} \rng(e_i)$ is also a homeomorphism.
			We put $Φ := \ImageMap{h}\restr{\H}$ and $\F := \rng(Φ)$.
			Clearly, $Φ\maps \H \to \F$ is a homeomorphism and $Φ(H) \homeo H$ for every $H ∈ \H$.
			
			For every $H ∈ \H$ and $i ∈ I$ we have $e_i\im{H ∩ V_i} ∈ \U_i$.
			This is because $e_i\im{H ∩ V_i} ⊆ Q_i ⊆ U_i$ and $H ∩ V_i ∈ ⋂_{j < s(i)} A_{i, j}^-$ so $e_i\im{H ∩ V_i} ∈ ⋂_{j < s(i)} U_{i, j}^-$.
			It follows that $Φ(H) = ⋃_{i ∈ I} e_i\im{H ∩ V_i} ∈ \U$, and so $\F ⊆ \U$.
		\end{proof}
	\end{proposition}
	
	Now we are ready to improve Proposition~\ref{thm:exact_shape_contains_space} from spaces to compact families of spaces.
	
	\begin{proposition} \label{thm:exact_shape_contains_family}
		Let $s\maps I \to ℕ_+$ be a finite function.
		For every strongly compactifiable class $\C ⊆ \O_s$ and every open set $\U ⊆ \Compacta{[0, 1]^ω}$ exactly of the shape $s$ there is a compact zero-dimensional family $\F ⊆ \U$ equivalent to $\C$.
		
		\begin{proof}
			By Theorem~\ref{thm:strongly_compactifiable_characterization} there is a closed zero-dimensional family $\H ⊆ \Compacta{[0, 1]^ω}$ equivalent to $\C$.
			For every $H ∈ \H$ let $\V_H ⊆ \Compacta{[0, 1]^ω}$ be a neighborhood of $F$ of the shape $s$ (Observation~\ref{thm:special_open_class}).
			The collection $\set{\V_H: H ∈ \H}$ is an open cover of $\H$.
			Since $\H$ is compact and zero-dimensional, there is a finite clopen decomposition $\set{\H_k: k < n}$ of $\H$ and a finite subcover $\set{\V_k: k < n} ⊆ \set{\V_H: H ∈ \H}$ such that $\H_k ⊆ \V_k$ for every $k < n$.
			
			By Proposition~\ref{thm:exact_shape_transport} for every $k < n$ there is homeomorphism $Φ_k\maps \H_k \to \F_k ⊆ \U$ such that $\H_k$ is equivalent to $\F_k$.
			Clearly, $\F := ⋃_{k < n} \F_k ⊆ \U$ is a compact zero-dimensional family equivalent to $\C$.
		\end{proof}
	\end{proposition}
	
	\begin{corollary}
		For every strongly compactifiable class of infinite compacta $\C$ and $ε > 0$ there is a closed zero-dimensional family $\F ⊆ \Compacta{[0 , 1]^ω}$ equivalent to $\C$ such that every space $F ∈ \F$ is $ε$-dense in $[0, 1]^ω$.
		
		\begin{proof}
			Let $A ⊆ [0, 1]^ω$ be a finite $2ε/3$-dense $2ε/3$-separated set and let $\U := ⋂_{x ∈ A} B(x, ε/3)^-$.
			The balls $B(x, ε/3)$ are pairwise disjoint, and so the open set $\U$ is exactly of the shape $s := \tuple{\card{A}}$.
			We have $\C ⊆ \O_s$ since all members of $\C$ are infinite and $\O_s$ is the class of all metrizable compacta with at least $\card{A}$ points.
			By Proposition~\ref{thm:exact_shape_contains_family} there is a closed zero-dimensional family $\F ⊆ \U$ equivalent to $\C$.
			For every $F ∈ \F$ and $x ∈ A$ we have $F ∩ B(x, ε/3) ≠ ∅$, and hence $F$ is $ε$-dense.
		\end{proof}
	\end{corollary}
	
	\begin{theorem} \label{thm:F_sigma_classes}
		Every countable union of strongly compactifiable classes is strongly compactifiable, i.e.\ every $F_σ$ subset of $\Compacta{[0, 1]^ω}$ is strongly compactifiable and equivalent to a closed subset of $\Compacta{[0, 1]^ω}$.
		
		\begin{proof}
			Let $\C_n$, $n ∈ ω$, be strongly compactifiable classes and let $\C = ⋃_{n ∈ ω} \C_n$.
			For every $n ∈ ω$ there is a compact zero-dimensional family $\H_n ⊆ \Compacta{[0, 1]^ω}$ equivalent to $\C_n$ (Theorem~\ref{thm:strongly_compactifiable_characterization}).
			The set of minimal types $R := A(\set{t(X): X ∈ \C})$ is finite as any antichain in $T ∪ \set{∞}$.
			For every $t ∈ R$ let us fix a space $F_{t, ∞} ∈ \Compacta{[0, 1]^ω}$ such that $F_{t, ∞} ∈ \homeocopies{\C}$ and $t(F_{t, ∞}) = t$.
			Every space $F_{t, ∞}$ has a countable decreasing neighborhood base $\set{\B_{t, n}: n ∈ ω}$ such that every $\B_{t, n}$ is exactly of the shape $s_{t, n}$ for some $s_{t, n} ∈ S_t$ (Lemma~\ref{thm:exact_shape_basis}).
			
			For every $n ∈ ω$ the family $\set{\O_{s_{t, n}}: t ∈ R}$ covers the compact zero-dimensional family $\H_n$ by Proposition~\ref{thm:principal_upper_class}, and so there is a clopen decomposition $\set{\H_{t, n}: t ∈ R}$ of $\H_n$ such that $\H_{t, n} ⊆ \O_{s_{t, n}}$ for every $t ∈ R$.
			By Proposition~\ref{thm:exact_shape_contains_family} there is a compact family $\F_{t, n} ⊆ \B_{t, n}$ equivalent to $\H_{t, n}$ for every $t ∈ R$.
			We put $\F_t := ⋃_{n ∈ ω} \F_{t, n} ∪ \set{F_{t, ∞}}$ and $\F := ⋃_{t ∈ R} \F_t$.
			Every family $\F_t$ is closed since the families $\F_{t, n}$ are closed and $⋂_{n ∈ ω} \clo{⋃_{m ≥ n} \F_{t, m}} ⊆ ⋂_{n ∈ ω} \clo{\B_{t, n}} = \set{F_{t, ∞}}$.
			The theorem follows since $\C = ⋃_{n ∈ ω} \C_n \homeo ⋃_{n ∈ ω} \H_n = ⋃_{t ∈ R, n ∈ ω} \H_{t, n} \homeo ⋃_{t ∈ R} \F_t = \F$.
		\end{proof}
	\end{theorem}
	
	\begin{corollary}
		Every $F_σ$ subset of $\Continua{[0, 1]^ω}$ is strongly compactifiable and equivalent to a closed subset of $\Continua{[0, 1]^ω}$.
	\end{corollary}

	Theorem~\ref{thm:F_sigma_classes} together with Theorem~\ref{thm:analytic_G_delta} and \ref{thm:open_classes} completes the picture of Borel complexity up to the equivalence – see Figure~\ref{fig:complexities}.
	The complexities reduce to four nontrivial groups of classes – clopen classes, open classes, strongly compactifiable classes, and strongly Polishable classes.

	
	\definecolor{R}{rgb}{0, 0, 0}
	\definecolor{G}{rgb}{0, 0, 0}
	\definecolor{B}{rgb}{0, 0, 0}
	\definecolor{Y}{rgb}{0, 0, 0}
	
	\pgfdeclaredecoration{allow_raise}{initial}{
		\state{initial}[width=\pgfdecoratedpathlength-1sp]{
			 \pgfmoveto{\pgfpointorigin}
		}
		\state{final}{
			 \pgflineto{\pgfpointorigin}
		}
	}

	\begin{figure}[!ht]
		\centering
		
		\begin{tikzpicture}[
				x = {(8em, 0)}, 
				y = {(0, 6.5em)},
				multiline/.style = {align=center}, 
				parallel/.style = {
					decoration = {allow_raise, raise=#1},
					decorate,
				},
				strong/.style = {parallel=-0.2em},
				weak/.style = {dashed, parallel=0.2em},
			]
			
			\node at (3.2, -0.5) (coanalytic) {$Π^1_1$ (co-analytic)};
			\node[outer sep=0.3em] at (3.7, 0) (last_dots) {$\cdots$};
			
			\begin{scope}[Y, multiline]
				\node at (0, 0) (clopen) {$Δ^0_1$ (clopen)};
				
				\node at (0, 1.2) (clopen_classes) {$∅, \set{∅}, \AllCompacta \setminus \set{∅}, \AllCompacta$};
				
				\graph{
					(clopen) <->[dashed] (clopen_classes)
				};
			\end{scope}
			
			\begin{scope}[R]
				\node at (0.5, 0.5) (open) {$Σ^0_1$ (open)};
				
				\node at (1.2, 1.2) (open_classes) {$\O_R\colon R ∈ \R$};
				
				\graph{
					(open) <->[dashed] (open_classes)
				};
			\end{scope}
			
			\begin{scope}[B, multiline]
				\node at (0.5, -0.5) (closed) {$Π^0_1$ (closed)};
				\node at (1, 0) (delta_2) {$Δ^0_2$};
				\node at (1.5, 0.5) (F_sigma) {$Σ^0_2$ ($F_σ$)};
				
				\node at (0.5, -1.2) (compactifiable) {strongly  compactifiable};
				
				\graph{
					(closed) ->[parallel=-0.2em] (delta_2) ->[parallel=-0.2em] (F_sigma),
					(closed) <-[dashed, parallel=0.2em] (delta_2) <-[dashed, parallel=0.2em] (F_sigma),
					(closed) <->[dashed] (compactifiable),
				};
			\end{scope}
			
			\begin{scope}[G, multiline]
				\node at (1.5, -0.5) (G_delta) {$Π^0_2$ ($G_δ$)};
				\node[outer sep=0.3em] at (2, 0) (dots) {$\cdots$};
				\node at (2.7, 0) (Borel) {$Δ^1_1$ (Borel)};
				\node at (3.2, 0.5) (analytic) {$Σ^1_1$ (analytic)};
				
				\node at (2.2, -1.2) (Polishable) {strongly  Polishable};
				
				\graph{
					(Borel) ->[strong] (analytic),
					(Borel) <-[weak] (analytic),
					(G_delta) <->[dashed] (Polishable),
					
					(G_delta) ->[strong] (dots) ->[strong] (Borel),
					(G_delta) <-[weak] (dots),
					(Borel) ->[strong, dashed] (dots), 
				};
			\end{scope}
			
			\graph{
				(closed) <- (clopen) -> (open) -> (delta_2) -> (G_delta),
				(F_sigma) -> (dots),
				(Borel) -> (coanalytic) -> (last_dots) <- (analytic),
			};
		\end{tikzpicture}
		
		\caption[Complexities and corresponding classes.]{Complexities and corresponding classes.
			“\tikz[baseline=-0.6ex]{\draw[->] (0, 0) -- (0.6, 0);}” denotes implication, 
			“\tikz[baseline=-0.6ex]{\draw[->, dashed] (0, 0) -- (0.6, 0);}” denotes implication up to the equivalence.
			}
		\label{fig:complexities}
	\end{figure}

	It is easy to see that there are open classes which are not clopen and that there are strongly compactifiable classes that are nor open. Also, there are classes which are not strongly Polishable. Nevertheless, the following remains open.
	
	\begin{question}
		Is there an analytic subset of $\Compacta{[0, 1]^ω}$ that is not equivalent to a closed subset?
		In other words, is there a class of metrizable compacta that is strongly Polishable, but not strongly compactifiable?
		One candidate is the class of all Peano continua \compactifiablePeano.
	\end{question}

\section{Saturated and type-saturated classes}
	
	We have defined saturated classes and saturated families.
	In general, on any set or class $X$ endowed with an equivalence we may consider its saturated subsets or subclasses – $A ⊆ X$ is saturated if it is the union of some equivalence classes, i.e. if it is closed under equivalent elements.
	So our saturated classes are saturated with respect to the equivalence of topological spaces where two spaces are equivalent if they are homeomorphic, and our saturated families are saturated with respect to the same equivalence but restricted to $\Compacta{[0, 1]^ω}$.
	
	\begin{definition}
		We say that a class of metrizable compacta $\C$ is \emph{type-saturated} if is it saturated with respect to the equivalence induced by the type function $t\maps \AllCompacta \to T ∪ \set{∞}$, i.e. $X, Y ∈ \AllCompacta$ are equivalent if $t(X) = t(Y)$.
		That means type-saturated classes are the unions $⋃_{t ∈ R} \T_t$ for $R ⊆ T ∪ \set{∞}$ where $\T_t$ for $t ∈ T ∪ \set{∞}$ denotes the \emph{principal type-saturated class} $\set{K ∈ \AllCompacta: t(K) = t}$.
		For a set of types $R ⊆ T ∪ \set{∞}$ we denote the type-saturated class $\set{K ∈ \AllCompacta: t(K) ∈ R} = ⋃_{t ∈ R} \T_t$ by $\T_R$.
		
		Clearly, every type-saturated class is saturated.
	\end{definition}
	
	\begin{remark}
		For every saturated class $\C$ of metrizable compacta we have $\homeocopies{(\C ∩ \Compacta{[0, 1]^ω})} = \C$,
		and for every saturated family $\F ⊆ \Compacta{[0, 1]^ω}$ we have $\homeocopies{\F} ∩ \Compacta{[0, 1]^ω} = \F$.
		This gives us a canonical identification between saturated classes and saturated families of metrizable compacta.
		Therefore, we may lift topological properties of saturated families to the corresponding saturated classes, e.g. we may say “closed class” or “open class” in the sense that the corresponding saturated family is closed or open.
		Note that this usage of “open class” is consistent with Definition~\ref{def:open_class}.
		This also includes the type-saturated classes, so for example “$\T_∞$ is $G_δ$” means that the corresponding family $\T_∞ ∩ \Compacta{[0, 1]^ω}$ is $G_δ$ in $\Compacta{[0, 1]^ω}$.
		On the other hand, we have defined only type-saturated classes, but this correspondence allows us to talk about type-saturated families without an explicit definition.
	\end{remark}
	
	\begin{observation} \label{thm:open_type-saturated}
		By Theorem~\ref{thm:open_classes} every open family $\U ⊆ \Compacta{[0, 1]^ω}$ is equivalent to some open class $\O_R$, which is by definition type-saturated. Hence, $\homeocopies{\U} = \O_R$.
		It follows that the saturation of an open family is still an open family, and that every saturated open or closed family is type-saturated.
		In particular, for a class $\C$ of metrizable compacta, $\homeocopies{\C} ∩ \Compacta{[0, 1]^ω}$ is closed if and only if $\homeocopies{\C} = \AllCompacta \setminus \O_R$ for some $R ∈ \R$.
		
		By Proposition~\ref{thm:clopen_classes} the situation with clopen families is even simpler – they just are type-saturated.
	\end{observation}
	
	The following corollary summarizes which complexities are preserved by saturation.
	
	\begin{corollary}
		If a family $\F ⊆ \Compacta{[0, 1]^ω}$ is clopen, open, or analytic, then so is the corresponding saturated family $\homeocopies{\F} ∩ \Compacta{[0, 1]^ω}$.
		On the other hand, there is a closed family $\F$ such that the corresponding saturated family is not Borel.
		
		\begin{proof}
			For clopen and open families, this follows Observation~\ref{thm:open_type-saturated}.
			The saturation of an analytic family is analytic by \compactifiableAnalytic{} and Theorem~\ref{thm:strongly_Polishable_characterization}.
			
			The class of all uncountable metrizable compacta is analytically complete \cite[Theorem~27.5]{Kechris}, but yet strongly compactifiable \compactifiableExample, and so equivalent to a closed family $\F$.
		\end{proof}
	\end{corollary}
	
	Let us make some remarks on the complexity of the saturation of a singleton family.
	So let $X$ be a metrizable compactum and let $\F$ be the corresponding saturated family $\homeocopies{\set{X}} ∩ \Compacta{[0, 1]^ω}$.
	$\F$ is always Borel \cite[Theorem~2]{RN_65}, but besides that it can be arbitrarily complex \cite[Fact~3.12]{Marcone}.
	Section~3.5 of \cite{Marcone} also gives us some examples:
	\begin{itemize}
		\item If $X$ is a graph or a dendrite with finitely many branching points, then $\F$ is $F_{σδ}$-complete \cite{CDM_05}.
		\item If $X$ is the pseudo-arc, then $\F$ is $G_δ$-complete since it is $G_δ$ and dense \cite{Bing_51}.
		\item If $X$ is the Sierpiński universal curve or the Menger universal curve, then $\F$ is $F_{σδ}$-complete \cite{Krupski_02}.
	\end{itemize}
	
	\begin{observation} \label{thm:singleton_saturation}
		It follows from Proposition~\ref{thm:neighborhood_types} that $\clo{\F} = \set{K ∈ \Compacta{[0, 1]^ω}: t(K) ≤ t(X)}$.
		Therefore, $\F$ is closed if and only if $X$ is degenerate.
		$\F$ is dense in nonempty compacta if and only if $t(X) = ∞$, i.e. if $X$ has infinitely many components.
		$\F$ is dense in nonempty continua if and only if $X$ is a nondegenerate continuum.
	\end{observation}

	In the last part we shall look at the type-saturated classes in more detail.
	We say that a type-saturated class $\T_R$ is \emph{lower} or \emph{upper} if the corresponding set $R$ is lower or upper in the ordered set $T ∪ \set{∞}$.
	Observe that every open type-saturated class is upper, and every closed type-saturated class is lower.
	
	Also recall that a subset of a topological space is called \emph{locally closed} if it is the intersection of an open set and a closed set.
	
	\begin{observation} \label{thm:principal_type-saturated}
		The class $\T_{0, 0}$ is clopen, $\T_{1, 0}$ is closed, $\T_∞$ is $G_δ$, and $\T_t$ is locally closed 
		for every other $t ∈ T ∪ \set{∞}$.
		No principal type-saturated class has a lower complexity than stated.
		
		\begin{proof}
			We already know that $\T_{0, 0} = \U_{0, 0} = \set{∅}$ is (with its complement) the only nontrivial clopen class (Proposition~\ref{thm:clopen_classes}).
			We have $\T_{1, 0} = \AllCompacta \setminus (\O_∅ ∪ \O_{\tuple{2}})$, so it is closed.
			We already know that $\T_∞ = \U_∞$ is $G_δ$ (Proposition~\ref{thm:principal_upper_class}) and dense (Observation~\ref{thm:singleton_saturation}), and so it is comeager.
			Since finite spaces are dense, $\T_∞$ has empty interior. So if it was $F_σ$, it would be also meager.
			For $t = \tuple{m, n} ∈ T_+$ we put $t' := \tuple{m, n + 1}$ if $m > n$ and $\tuple{m + 1, 0}$ otherwise.
			Let $\V := \U_t ∪ \U_{\tuple{m + 1, 0}}$ and $\V' := \U_{t'} ∪ \U_{\tuple{m + 1, 0}}$.
			Both classes $\V$ and $\V'$ are open and $\T_t = \V \setminus \V'$, so $\T_t$ is locally closed.
			$\T_t$ for $t ∉ \set{\tuple{0, 0}, \tuple{1, 0}, ∞}$ is neither open nor closed since it is neither upper nor lower.
		\end{proof}
	\end{observation}
	
	\begin{corollary}
		Let $R ⊆ T ∪ \set{∞}$.
		If $∞ ∉ R$, then $\T_R$ is $F_σ$.
		Otherwise, $\T_R$ is $G_δ$.
		
		\begin{proof}
			We have $\T_R = ⋃_{t ∈ R} \T_t$, and if $∞ ∉ R$, then each such $\T_t$ is $F_σ$.
			If $∞ ∈ R$, then the complementing type-saturated class is $F_σ$ by the previous claim.
		\end{proof}
	\end{corollary}
	
	\begin{remark}
		Even though the class $\U_∞ = \T_∞$ of all metrizable compacta with infinitely many components is not $F_σ$, it is strongly compactifiable \compactifiableInftyComponents.
		It follows that every type-saturated class $\T_R$, $R ⊆ T ∪ \set{∞}$, is strongly compactifiable since it is either $\T_{R \setminus \set{∞}}$ or $\T_{R \setminus \set{∞}} ∪ \T_∞$, and $\T_{R \setminus \set{∞}}$ is $F_σ$ by the previous corollary.
	\end{remark}
	
	\begin{remark}
		In the previous corollary we used the fact that every open saturated family is $F_σ$. But that does not mean it is the countable union of saturated closed families.
		Saturated closed families are type-saturated (Observation~\ref{thm:open_type-saturated}), so every union of them is also type-saturated.
		On the other hand, there are $F_σ$ or $G_δ$ saturated families that are not type-saturated (see the examples before Observation~\ref{thm:singleton_saturation}).
	\end{remark}
	
	\begin{observation}
		Let us consider the quotient $q_{\homeo}\maps \Compacta{[0, 1]^ω} \to \Compacta{[0, 1]^ω}/{\homeo}$, so open subsets of $\Compacta{[0, 1]^ω}/{\homeo}$ correspond to saturated open families.
		In general, subsets of $\Compacta{[0, 1]^ω}/{\homeo}$ correspond to saturated families, and for example $F_σ$ subsets of $\Compacta{[0, 1]^ω}/{\homeo}$ correspond to countable unions of saturated closed families.
		Since by the proof of Observation~\ref{thm:principal_type-saturated} every principal type-saturated class is obtained as a Borel combination of open type-saturated classes, we have that type-saturated classes correspond exactly to Borel subsets of $\Compacta{[0, 1]^ω}/{\homeo}$.
		
		It is not true that open subsets of $\Compacta{[0, 1]^ω}/{\homeo}$ are $F_σ$. This space is not metrizable. In fact, it is not even $T_0$.
		Two points of $\Compacta{[0, 1]^ω}/{\homeo}$ represented by spaces $X, Y ∈ \Compacta{[0, 1]^ω}$ are indistinguishable if and only if $t(X) = t(Y)$, so we may consider the Kolmogorov quotient $q_{T_0}\maps \Compacta{[0, 1]^ω}/{\homeo} \to T ∪ \set{∞}$.
		In fact, the composition quotient map $q_{T_0} ∘ q_{\homeo}$ is just the type function $t\maps \Compacta{[0, 1]^ω} \to T ∪ \set{∞}$.
		This endows the set of all types $T ∪ \set{∞}$ with the topology where $R ⊆ T ∪ \set{∞}$ is open if and only if it is upper and nice.
		
		It is also easy to directly see that these sets form a topology.
		Upper sets are stable under arbitrary unions and intersections, and nice sets are stable under arbitrary unions.
		Moreover, nice upper sets are stable under finite intersections: if $R_1 ∩ R_2 ∩ T_+ ≠ ∅$, then since $R_1$ and $R_2$ are nice, there are some $m_1, m_2 > 0$ such that $\tuple{m_1, 0} ∈ R_1$ and $\tuple{m_2, 0} ∈ R_1$. Since $R_1$ and $R_2$ are upper, we have $\max\set{\tuple{m_1, 0}, \tuple{m_2, 0}} ∈ R_1 ∩ R_2$.
	\end{observation}
	
	\begin{observation}
		The proof of Observation~\ref{thm:principal_type-saturated} in fact works in $T ∪ \set{∞}$, i.e. $\set{\tuple{0, 0}}$ is clopen, $\set{\tuple{1, 0}}$ is closed, $\set{∞}$ is $G_δ$, and $\set{t}$ is locally closed for every other $t ∈ T ∪ \set{∞}$.
		Also, no singleton has a lower complexity than stated.
		
		Here we have to be more careful since open sets are not necessarily $F_σ$.
		Instead of $F_σ$ we should consider the complexity $Σ^0_2$ – the countable unions of members of $Π^0_1$.
		Instead of starting just with open sets and closed sets, we let $Π^0_1 = Σ^0_1$ be the algebra generated by open sets and closed sets.
		Members of the algebra are called \emph{constructible sets}, and they are finite unions of locally closed sets.
		
		So let us show that $\set{∞}$ is not $Σ^0_2$. Since our set is a singleton, it would mean $\set{∞}$ is locally closed.
		If $\set{∞}$ was locally closed in $T ∪ \set{∞}$, we would have $\set{∞} = \clo{\set{∞}} ∩ U = (T_+ ∪ \set{∞}) ∩ U = U$ for some open set $U ⊆ T_+ ∪ \set{∞}$. So $\set{∞}$ would be open, which it is not since it is not nice.
		
		Also, for $t ≠ \tuple{0, 0}, \tuple{1, 0}, ∞$ the singleton $\set{t}$ is neither in any class $F_σ$, $F_{σδ}$, $F_{σδσ}$, … since they consist only of lower sets, nor in any class $G_δ$, $G_{δσ}$, $G_{δσδ}$, … since they consist only of upper sets.
	\end{observation}

	\begin{acknowledgements}
		The author is grateful to his supervisor Benjamin Vejnar for discussions and proofreading and also to Włodzimierz J.\ Charatonik for the suggestion that strongly compactifiable classes could be stable under countable unions and that this could be proved by approximating a limit object in the hyperspace.
		
		The author was supported by the grant projects GAUK 970217 and SVV-2017-260456 of Charles University.
	\end{acknowledgements}
	
	\linespread{1}\selectfont
	
	\bibliographystyle{mysiam}
	\bibliography{references}
	
\end{document}